\documentclass[11pt,reqno]{amsart}
\usepackage{amssymb}
\newcommand{\la}{\lambda}
\newcommand{\al}{\alpha}

\newcommand{\p}{\partial}
\newcommand{\I}{\infty}

\newcommand{\ti}{\tilde}
\newcommand{\R}{\mathbb{R}}

\newcommand{\Z}{\mathbb{Z}}
\newcommand{\F}{\mathcal{F}}

\numberwithin{equation}{section}
\newtheorem{thm}{Theorem}

\newtheorem{prop}[thm]{Proposition}

\theoremstyle{remark}
\newtheorem{rem}{Remark}

\newcommand{\lec}{\lesssim}

\newcommand{\EQ}[1]{\begin{equation} \begin{split} #1 \end{split} \end{equation}}
\newcommand{\Del}[1]{}

\newcommand{\pq}{\quad}

\newcommand{\de}{\delta}

\newcommand{\x}{\xi}

\newcommand{\na}{\nabla}

\newcommand{\supp}{\operatorname{supp}}

\begin{document}
\title{Sharp global regularity for the $2+1$-dimensional equivariant Faddeev model}

\author{Dan-Andrei Geba, Kenji Nakanishi, and Xiang Zhang}

\address{Department of Mathematics, University of Rochester, Rochester, NY 14627, U.S.A.}
\email{dangeba@math.rochester.edu}
\address{Department of Mathematics, Kyoto University, Kyoto 606-8502, Japan}
\email{n-kenji@math.kyoto-u.ac.jp}
\address{Department of Mathematics, University of Rochester, Rochester, NY 14627, U.S.A.}
\email{zhang@math.rochester.edu}
\date{}

\begin{abstract}
The aim of this article is to prove that for the $2+1$-dimensional equivariant Faddeev model, which is a quasilinear generalization of the corresponding nonlinear $\sigma$ model, small initial data in critical Besov spaces evolve into global solutions which scatter. 
\end{abstract}


\maketitle

\section{Introduction}

\subsection{Physical background and previous results} 
A classical field theory which models elementary heavy particles by topological solitons was introduced by Faddeev in \cite{F75, F76}. One of the remarkable features of the Faddeev model is that it admits knotted solitons. It is described by the action
\begin{equation}
S\,=\,  \int_{\R^{3+1}}\,\frac{1}{2}\partial_\mu {\bf n}\cdot \partial^\mu {\bf n}\, + \,\frac{1}{4}(\partial_\mu {\bf n}\wedge\partial_\nu {\bf n})\cdot(\partial^\mu {\bf n}\wedge\partial^\nu {\bf n})\,dg,
\label{fdv}
\end{equation}
where $v_1\wedge v_2$ denotes the cross product of the vectors $v_1$ and $v_2$ in $\mathbb{R}^3$ and ${\bf n}: \R^{3+1}\to \mathbb{S}^2$ are maps from the Minkowski spacetime, with $g= \text{diag}(-1,1,1,1)$, into the unit sphere of $\mathbb{R}^{3}$. Formal critical points for this theory can be described variationally by  
\begin{equation}
{\bf n}\wedge\partial_\mu\partial^\mu {\bf n}+(\partial_\mu[{\bf n}\cdot(\partial^\mu {\bf n}\wedge\partial^\nu {\bf n})])\partial_\nu {\bf n}=0,
\label{fsys}
\end{equation}
which is a system of quasilinear wave equations. Natural extensions of this model can be obtained by switching   the domain of the maps ${\bf n}$ from $\R^{3+1}$ to $\R^{n+1}$, endowed also with the Minkowski metric. 

The Faddeev model is profoundly related with the first classical theory which models particles by topological solitons. This is the Skyrme model \cite{S1, S2, S3}, whose action is given by
\begin{equation}
\label{sk}
S\,=\,\int_{\mathbb{R}^{3+1}} \frac 12\,\langle\partial^\mu\phi\,,\, \partial_\mu\phi \rangle_h\,+\frac{\alpha^2}{4}\,\left( \langle\partial^\mu\phi\,,\, \partial_\mu\phi \rangle_h^2-\langle\partial^\mu\phi\,,\, \partial^\nu\phi \rangle_h \langle\partial_\mu\phi\,,\, \partial_\nu\phi \rangle_h \right)\,dg,
\end{equation}
where $\alpha$ is a constant having the dimension of length and $\phi: \R^{3+1}\to \mathbb{S}^3$ are maps from the Minkowski spacetime into the unit sphere of $\mathbb{R}^{4}$. If one restricts the image of $\phi$ to be the equatorial 2-sphere of $\mathbb{S}^3$ (identified in this case with $\mathbb{S}^2$) by prescribing 
\begin{equation}
\phi=(u, {\bf n})=(\pi/2, {\bf n}),
\end{equation} 
where the metric on $\mathbb{S}^3$ is 
\begin{equation}
 h=du^2+\sin^2 u\,d{\bf n}^2,\qquad 0 \leq u\leq \pi, \qquad {\bf n}\in\mathbb{S}^2,
\end{equation}
then the Skyrme action \eqref{sk} reduces to the Faddeev one \eqref{fdv}. Moreover, the $2+1$-dimensional Skyrme model, which is given by 
\begin{equation}
\label{sk2}
S\,=\,\int_{\mathbb{R}^{2+1}} \frac 12\,\langle\partial^\mu\phi\,,\, \partial_\mu\phi \rangle_h\,+\frac{\alpha^2}{4}\,\left( \langle\partial^\mu\phi\,,\, \partial_\mu\phi \rangle_h^2-\langle\partial^\mu\phi\,,\, \partial^\nu\phi \rangle_h \langle\partial_\mu\phi\,,\, \partial_\nu\phi \rangle_h \right)\,+V(\phi) \,dg,
\end{equation}
becomes the $2+1$-dimensional Faddeev theory in the particular case when the potential term vanishes.

Initial studies on the Faddeev model focused on its static properties, with numerical simulations for various topological solitons being performed by Faddeev-Niemi \cite{FN97} and Battye-Sutcliffe \cite{BS98, BS99}, while Vakulenko-Kapitanski \cite{VK79} and Lin-Yang \cite{LY104, LY04} investigated the associated topologically-constrained energy-minimization problem. To our knowledge, the only result for the time evolution of the Faddeev model is due to Lei-Lin-Zhou \cite{LLZ}, who proved that the $2+1$-dimensional system \eqref{fsys} is globally well-posed (GWP) for smooth, compactly-supported initial data with small $H^{11}(\R^2)$ norm.

\subsection{Motivation and formulation of the problem} 
In this article, we study the global regularity for the $2+1$-dimensional equivariant Faddeev model. Our motivation is threefold. First, we would like to lower the regularity needed for GWP in Lei-Lin-Zhou's result. Secondly, the $2+1$-dimensional Faddeev model is a quasilinear generalization of the energy-critical wave maps system, which has long been a hot topic in the field of hyperbolic equations. Finally, our choice for the equivariant assumption has to do with the fact that this is the most natural rotational symmetry one can associate to maps ${\bf n}: \R^{2+1}\to \mathbb{S}^2$. 

Thus, our focus is on equivariant maps ${\bf n}: (\mathbb{R}^{2+1},g)\to(\mathbb{S}^2,h)$, i.e.,
\begin{equation}
{\bf n}(t,r,\omega)\,=\,(u(t,r),\omega), \quad g= -dt^2+dr^2+r^2\,d\omega^2,\quad h=du^2+\sin^2u\,d\omega^2,
\end{equation}
which are critical points for the $2+1$-dimensional Faddeev model. Using this ansatz, the Euler-Lagrange system \eqref{fsys} simplifies to the following quasilinear wave equation, satisfied by the angular variable $u$:
\begin{equation}
\label{meq}
\left(1+\frac{\sin^2u}{r^2}\right) (u_{tt}-u_{rr}) - \left(1-\frac{\sin^2u}{r^2}\right) \frac{u_r}{r} + \frac{\sin 2u}{2r^2}\left(u_t^2-u_r^2+1\right)= 0.
\end{equation}
The a priori conserved energy associated to this equation is given by
\begin{equation}
E[u]\,=\,\int_0^\infty \left[\left(1 + \frac{\sin ^2 u}{r^2}\right)\frac{u_t^2+u_r^2}{2}+
\frac{\sin ^2 u}{2r^2} 
\right]\,r dr. \label{tote}
\end{equation}

If we choose to work with finite energy maps, then we would also need to assume that $\sin u(t,0)=\sin u(t,\infty)=0$, which can be satisfied by restricting our study to degree-0 equivariant maps, i.e.,
\begin{equation}
u(t,0)\,=\,u(t,\infty)\,=\,0.
\end{equation}
Using this extra hypothesis, we can show that smooth, finite energy solutions of \eqref{meq} are uniformly bounded  and 
\EQ{\label{li} \|u\|_{L^\infty_{t,x}} \leq C(E[u]),}
with $C(s)\to 0$ as $s\to 0$. This follows by introducing
\EQ{
I(z)=\int_0^z |\sin w|\,dw,}
which verifies
\EQ{
I(0)=0, \qquad |I(z)|>0\,(z\neq 0),\quad \text{and} \quad \lim_{|z|\to\infty}
|I(z)|=\infty,}
and then arguing as follows:
\begin{equation}
\aligned
|I(u(t,r))|\,&=\,\left|\int_0^r |\sin u(t,s)|\, u_r(t,s)\, ds\right|\\
&\lesssim \,\left(\int_0^r \, \frac{\sin ^2 u(t,s)}{s}\,ds \right)^\frac 12 \
\left(\int_0^r \,u^2_r(t,s) s\,ds\right)^\frac 12\,\lesssim\,E(u).
\endaligned
\label{bog}\end{equation}

The quasilinear equation is not scale-invariant. However, if we work in the small energy regime and we use the bound \eqref{li}, then we can formally write a scale-invariant approximation,
\EQ{\label{smeq} \left(1+\frac{u^2}{r^2}\right)(u_{tt}-u_{rr})-\left(1-\frac{u^2}{r^2}\right)\frac{u_r}{r}+\frac{u}{r^2} \left(u_t^2-u_r^2+1\right)=0,}  
for which, if $u$ is a solution, then so is
\EQ{
u_{\la}(t,r) = \la\, u\left(\frac t\la, \frac r\la\right).}
As
\EQ{
 \|u_\la(0)\|_{\dot H^{2}(\R^2)} = \|u(0)\|_{\dot H^{2}(\R^2)},}
we are naturally led to believe that an optimal result for the Cauchy problem associated to \eqref{meq} is a GWP for small data in $\dot H^{2}(\R^2)$-like spaces. On the other hand, the energy \eqref{tote} is at $\dot H^{3/2} \cap \dot H^1(\R^2)$-regularity level, which tells us that the equation is supercritical with respect to the energy.

Finally, performing the substitution $u=rv$, we can rewrite \eqref{meq} as a semilinear wave equation for $v$ in $\R^{4+1}$,  
\begin{equation}
\label{fdveq}
\Box v=h_1(r,u)\cdot v^3v_r+h_2(r,u)\cdot v^3+h_3(r,u)\cdot v^5+h_4(r,u)\cdot v(v_t^2-v_r^2),
\end{equation}
where $\Box=-\p^2_{t}+\p^2_{r}+\frac{3}{r}\,\p_r$ is the radial wave operator on $\R^{4+1}$ and the coefficients of the nonlinearities are given by
\begin{equation}
 \begin{cases}
 \label{Nonlinear coefs}
 h_1(r,u)=\frac{2\sin u(\sin u-u\cos u)}{u^3\, \Phi(r,u)}, \quad h_2(r,u)=\frac{\sin 2u-2u}{2u^3\, \Phi(r,u)},\\
 h_3(r,u)=\frac{\sin u(\sin u-u\cos u)}{u^4\, \Phi(r,u)}, \quad \ \, h_4(r,u)=\frac{\sin 2u}{2u\, \Phi(r,u)},
 \end{cases}
\end{equation}
with
\begin{equation}
\Phi(r,u)\,=\,1+\frac{\sin^2 u}{r^2}.
\end{equation}

Using the notation
\begin{equation}
\tilde{h}_i(u)\,=\,h_i(r,u)\,\Phi(r,u),\qquad 1\leq i\leq 4,
\end{equation}
direct computations based on Maclaurin series yield the following decay estimates:
\begin{equation}
 \begin{cases}
 |\p_{u}^j \tilde{h}_1(u)|\,\lesssim\, \langle u\rangle^{-2}, \quad |\p_{u}^j \tilde{h}_3(u)|\,\lesssim\, \langle u\rangle^{-3}, \quad |\p_{u}^j \tilde{h}_4(u)|\,\lesssim\, \langle u\rangle^{-1},\\
 |\tilde{h}_2(u)|\,\lesssim\, \langle u\rangle^{-2},\qquad |\p_{u}^{1+j} \tilde{h}_2(u)|\,\lesssim\, \langle u\rangle^{-3}, 
\end{cases}
\label{analytic}
\end{equation}
for all integers $j\geq 0$, with $\langle u\rangle=(1+u^2)^{\frac 12}$.
 
\subsection{Main result} 
We can now state our main result, which validates the previous scaling heuristics.

\begin{thm}
There exists $\de>0$ such that for any radial intial data $(v(0,r),\p_t v(0,r))$ decaying as $r\to\I$ and satisfying 
\EQ{
 \|\p v(0,\cdot)\|_{\dot B^1_{2,1}\cap \dot B^0_{2,1}(\R^4)}\le\de,}
the equation \eqref{fdveq} admits a unique global solution $v$ verifying 
\EQ{\label{reg}
 \p v \in C(\R;\dot B^1_{2,1}\cap \dot B^0_{2,1}(\R^4))\cap L^2(\R;\dot B^{1/6}_{6,1}\cap\dot B^{- 5/6}_{6,1}(\R^4)).}
Moreover, for some $v_\pm$ solving the free wave equation, 
\EQ{
 \|\p(v-v_\pm)(t)\|_{\dot B^1_{2,1}\cap \dot B^0_{2,1}(\R^4)} \to 0 \pq \text{as}\pq t\to\pm\I\, .}
\label{mainth}
\end{thm}

\begin{rem} For $s\in \R$ and $1\leq p,q\leq \infty$, $\dot B^{s}_{p,q}(\R^n)$ represents the homogeneous Besov space which is defined using a classical Littlewood-Paley decomposition, i.e.,
\EQ{\label{B}
u\,=\,\sum_{\lambda\in 2^{\mathbb{Z}}} S_\lambda(\na) u\,=\,\sum_{\lambda\in 2^{\mathbb{Z}}} u_\la, \qquad \|u\|_{\dot B^{s}_{p,q}(\R^n)}\,=\,\left(\sum_{\lambda\in 2^{\mathbb{Z}}} \left(\lambda^{s}\,\|u_\la\|_{L^p(\R^n)}\right)^q\right)^{1/q},}
where $S_\la(\na)=\F^{-1}\chi(\la^{-1}\x)\F$ is the Fourier multiplier on $\R^n$ given by a fixed radial function $\chi\in C_0^\I(\R^n)$ satisfying 
\[
\supp\chi\subset\{1/2<|\x|<2\}, \qquad \sum_{\la\in 2^\Z}\chi(\la^{-1}\x)=1, \ (\forall)\x\not=0.\] 
The transition between the norms for $v$ and the ones for $u=rv$ should be done according to 
\begin{equation}
\|v(t)\|_{\dot{B}^s_{p,q}(\R^4)}\,\simeq\,\left \|r^{\frac 2p -1}\,u(t)\right\|_{\dot{B}^s_{p,q}(\R^2)}.
\end{equation}
\end{rem}

\begin{rem}
A similar result has been obtained by the third author, in his doctoral dissertation, for the action
\begin{equation}
S\,=\,  \int_{\R^{2+1}}\,\frac{1}{2}\partial_\mu {\bf n}\cdot \partial^\mu {\bf n}\, - \,\frac{1}{4}(\partial_\mu {\bf n}\wedge\partial_\nu {\bf n})\cdot(\partial^\mu {\bf n}\wedge\partial^\nu {\bf n})\,dg.
\label{fdvz}
\end{equation}
\end{rem}

\subsection{Method of proof}
For proving Theorem \ref{mainth} we use a contraction-based argument applied to the integral version of \eqref{fdveq},
\EQ{
v\,=\,S(v_0,v_1)\,+\,\Box^{-1}(N(v)),}
where $S=S(v_0,v_1)$ is the homogeneous solution operator,
\EQ{
\Box\, S\,=\,0, \qquad S(0)\,=\,v_0, \qquad S_t(0)\,=\,v_1,}
$\Box^{-1}=\Box^{-1}F$ is the inhomogeneous solution operator,
\EQ{
\Box\,\left(\Box^{-1}F\right)\,=\,F, \qquad \left(\Box^{-1}F\right)(0)\,=\,\left(\Box^{-1}F\right)_t(0)\,=\,0,}
and 
\EQ{N(v)=h_1(r,u)\cdot v^3v_r+h_2(r,u)\cdot v^3+h_3(r,u)\cdot v^5+h_4(r,u)\cdot v(v_t^2-v_r^2).}
This reduces the proof to finding a data space $D$ and a solution space $X$ for which the following estimates are true:
\begin{align}
&\|S(v_0,v_1)\|_X\,\lesssim\,\|(v_0,v_1)\|_D, \qquad \qquad \quad \,\|w\|_{L_t^\infty\,D}\,\lesssim\,\|w\|_X, \label{hom}\\
&\|\Box^{-1}(N(w_1)\,-\,N(w_2))\|_X\,\lesssim\,(\|w_1\|_X+\|w_2\|_X)\|w_1-w_2\|_X \label{lip},
\end{align}
where the last bound holds for $\|v_1\|_X$ and $\|v_2\|_X$ sufficiently small.

In our case, the data space $D$ is described by
\begin{equation}
\|(v_0,v_1)\|_D\,=\,\|v_0\|_{\dot{B}^{2}_{2,1}\cap\dot{B}^{1}_{2,1} (\R^4)}+\|v_1\|_{\dot{B}^{1}_{2,1}\cap\dot{B}^{0}_{2,1}(\R^4)},
\label{dnorm}
\end{equation}
while the solution space $X$ is the $n=4$ version of the function space $Z$ introduced by Geba-Nakanishi-Rajeev \cite{GNR} in proving a similar result for the $3+1$-dimensional equivariant Skyrme model. The space $Z$ is based on the Besov hyperbolic spaces developed by Tataru in \cite{T} for the GWP of higher-dimensional wave maps.

The homogeneous bounds \eqref{hom} can be derived immediately from the properties of Tataru's spaces. The most involved part of the argument is the proof of the nonlinear estimate \eqref{lip}. Based on the analytic properties \eqref{analytic} for the nonlinear coefficients $\tilde{h}_i$ and the polynomial structure of $N(v)$, we will show that, for $\|v\|_X$ sufficiently small,
\begin{equation}
\label{nonlin}
\|\Box^{-1}(N(v))\|_{X}\,\lesssim\,\|v\|^3_X.
\end{equation}
Then it will be clear that \eqref{lip} is derived by the same argument and ingredients, for which the detail will be omitted. 

The next section is devoted to the proof of 
\begin{equation}
\label{nonlin*}
\aligned
\|\tilde{N}(v)\|_{L^1_t\left(\dot{B}^1_{2,1}\cap\dot{B}^0_{2,1}(\R^{4})\right)}\,\lesssim \,\|\partial v\|^3_{L^{\infty}_t\left(\dot{B}^1_{2,1}\cap \dot{B}^{0}_{2,1}(\R^4)\right)}+\|\partial v\|^3_{L^2_t\left(\dot{B}^{1/6}_{6,1}\cap\dot{B}^{-5/6}_{6,1}(\R^4)\right)},
\endaligned
\end{equation}
where 
\begin{equation}
\tilde{N}(v)=h_1(r,u)\cdot v^3v_r+h_2(r,u)\cdot v^3+h_3(r,u)\cdot v^5.
\end{equation} 
This is a stronger bound than the one we need,
\begin{equation}
\label{nonlin3*}
\left\|\Box^{-1}(\tilde{N}(v))\right\|_{X}\,\lesssim\,\|v\|^3_X,
\end{equation}
as we show in the final section that $\Box^{-1}$ maps $L^1_t\left(\dot{B}^1_{2,1}\cap\dot{B}^0_{2,1}(\R^{4})\right)$ into $X$ and $\|v\|_X^3$ dominates the right-hand side of \eqref{nonlin*}.

The final section also deals with the last nonlinear term, 
\begin{equation}
h_4(r,u)\cdot v(v_t^2-v_r^2),
\label{nf}
\end{equation}
for which we use its null-form structure,
\begin{equation}
\label{q}
v_t^2-v_r^2=-\Box \left(\frac{v^2}{2}\right)+v\,\Box\, v.
\end{equation}
In proving
\begin{equation}
\label{nf3}
\left\|\Box^{-1}(h_4(r,u)\cdot v(v_t^2-v_r^2))\right\|_{X}\,\lesssim\,\|v\|^3_X,
\end{equation}
it would be enough for the bilinear inequality
\begin{equation}
\|r\,u\,v\|_Z\,\lesssim\,\|u\|_Z\,\|v\|_Z
\end{equation}
to be true. This was proved in \cite{GNR} for spatial dimensions $n\geq 5$ and, unfortunately, that argument can not be adapted to our context (i.e., $n=4$). Instead, and this is the main novelty of our paper, we demonstrate the trilinear bound
\begin{equation}
\|r^2\,u\,v\,w\|_Z\,\lesssim\,\|u\|_Z\,\|v\|_Z\,\|w\|_Z,
\end{equation}
which allows us to infer \eqref{nf3}.

\section{The analysis for the cubic, quintic, and quartic nonlinearities}
As mentioned previously, the goal of this section is to prove the nonlinear estimate \eqref{nonlin*}.

\subsection{Algebra-type Besov estimates}
First, a direct argument based on paradifferential calculus and Bernstein and H{\"o}lder inequalities leads to 
\begin{equation}
\|f\, g\|_{\dot{B}^s_{p,1}(\R^n)}\lesssim \|f\|_{L^{p_1}(\R^n)}\, \|g\|_{\dot{B}^s_{q_1,1}(\R^n)}+\|g\|_{L^{p_2}(\R^n)}\, \|f\|_{\dot{B}^s_{q_2,1}(\R^n)},
\label{Bes-prod}\end{equation}
if 
\begin{equation}
\quad \frac{1}{p}=\frac{1}{p_1}+\frac{1}{q_1}=\frac{1}{p_2}+\frac{1}{q_2},\quad p, p_1, p_2, q_1, q_2\geq 1, \quad s>0.
\end{equation}

Secondly, using the previous product bound and the standard homogeneous Besov embedding
\EQ{\label{Bes2}
\dot{B}^{s_1}_{p_1,1}(\R^n)\subset \dot{B}^{s_2}_{p_2,1}(\R^n), \quad s_1-s_2=n\left(\frac{1}{p_1}-\frac{1}{p_2}\right), \ 1\leq p_1 \leq p_2, \ s_1, s_2 \in\R,}
we deduce
\begin{equation}
\|f\, g\|_{\dot{B}^0_{2,1}(\R^4)}\lesssim \|f\|_{L^4(\R^4)}\, \|g\|_{\dot{B}^1_{2,1}(\R^4)}+\|g\|_{L^4(\R^4)}\, \|f\|_{\dot{B}^1_{2,1}(\R^4)}.
\label{b021}\end{equation}

This is followed by:
\begin{prop}
The space $Y=\dot{B}^2_{2,1}\cap \dot{B}^{1}_{2,1}(\R^4)$ is an algebra, i.e.,
\begin{equation}
\|w_1\,w_2\|_Y\,\lesssim\,\|w_1\|_Y\,\|w_2\|_Y.
\label{y}
\end{equation}
Moreover, for radial functions $w_1=w_1(r)$ and $w_2=w_2(r)$, we have
\begin{equation}
\|r\,w_1\,w_2\|_Y\,\lesssim\,\|w_1\|_Y\,\|w_2\|_Y.
\label{ry}
\end{equation} \label{Yprop}
\end{prop} 
\begin{proof}
The algebra inequality is a direct consequence of \eqref{Bes-prod}, with $p_1=p_2=\infty$, and the Sobolev embedding
\begin{equation}
\dot{B}^{n/2}_{2,1}(\R^n)\subset L^{\infty}(\R^n).
\label{sob}
\end{equation} 
In demonstrating the radial bound \eqref{ry}, we use a radial Sobolev inequality proved in \cite{GNR},
\begin{equation}
\label{pq}
 \|r^{\alpha(1/p-1/q)}\phi_\la\|_{L^q(\R^n)} \lesssim \la^{(n-\alpha)(1/p-1/q)} \|\phi\|_{L^p(\R^n)},
\end{equation}
which holds for all radial functions $\phi=\phi(r) \in L^p(\R^n)$, $0\leq\alpha\leq n-1$, and $2\leq p\leq q$. For $n=4$, $\alpha=p=2$, and $q=\infty$, this reads as
\begin{equation}
\label{pq4}
 \|r\phi_\la\|_{L^\infty(\R^4)} \lesssim \la \|\phi_\la\|_{L^2(\R^4)}.
\end{equation}
Due to this estimate, we can infer 
\begin{equation*}
\aligned
&\|r\,w_1\,w_2\|_{\dot{B}^2_{2,1}}\\ &\lesssim \sum\limits_{\la}\,\sum\limits_{\mu\lesssim \la}\,\la^{2}\left(\|r\,w_{1,\la}\,w_{2,\mu}\|_{L^2}+\|r\,w_{1,\mu}\,w_{2,\la}\|_{L^2}\right)+ \sum\limits_{\la}\,\sum\limits_{\mu\gtrsim \la} \,\la^{2}\|r\,w_{1,\mu}\,w_{2,\mu}\|_{L^2}\\
&\lesssim \sum\limits_{\la}\,\sum\limits_{\mu\lesssim \la}\,\la^{2}\left(\|w_{1,\la}\|_{L^2}\,\|r\,w_{2,\mu}\|_{L^\infty}+\|r\,w_{1,\mu}\|_{L^\infty}\,\|w_{2,\la}\|_{L^2}\right)+ \sum\limits_{\mu} \,\mu^{2}\|r\,w_{1,\mu}\|_{L^\infty}\,\|w_{2,\mu}\|_{L^2}\\
&\lesssim \sum\limits_{\la}\,\sum\limits_{\mu\lesssim \la}\,\mu\,\la^{2}\left(\|w_{1,\la}\|_{L^2}\,\|w_{2,\mu}\|_{L^2}+\|w_{1,\mu}\|_{L^2}\,\|w_{2,\la}\|_{L^2}\right)+ \sum\limits_{\mu} \,\mu^{3}\|w_{1,\mu}\|_{L^2}\,\|w_{2,\mu}\|_{L^2}\\
&\lesssim \|w_1\|_{\dot{B}^2_{2,1}}\, \|w_2\|_{\dot{B}^1_{2,1}} + \|w_1\|_{\dot{B}^1_{2,1}}\, \|w_2\|_{\dot{B}^2_{2,1}} \,\lesssim\,\|w_1\|_Y\,\|w_2\|_Y.
\endaligned
\end{equation*}
The argument for the $\dot{B}^1_{2,1}$ norm is identical. 
\end{proof}

Finally, we derive estimates related to the nonlinear coefficient $\Phi(r,u)$.
\begin{prop}
The following inequality
\begin{equation}
\left\|\left(\frac{\sin rv}{r}\right)^{2k}\right\|_Y \,\lesssim\,C^k \|v\|^{2k}_Y, \qquad (\forall) k\geq 1,
\label{ysin}\end{equation}
holds for all radial functions $v$ with $\|v\|_Y \leq 1$, where  $C>0$ is a positive constant independent of $v$.
\end{prop}
\begin{proof}
Based on the fact that $Y$ is an algebra, \eqref{ysin} follows if we prove 
\begin{equation}
\left\|\left(\frac{\sin rv}{r}\right)^{2}\right\|_Y \,\lesssim\,\|v\|^{2}_Y, 
\label{ysin2}\end{equation}
for $\|v\|_Y\leq 1$. However, this is immediate if we combine the Maclaurin series for $\sin^2 x$ with 
\begin{equation}
\|r^{2k}\, v^{2k+2}\|_Y \,\lesssim\,C^k \|v\|^{2k+2}_Y, \qquad (\forall) k\geq 0,
\label{yn}\end{equation}
which, in turn, can be deduced using \eqref{ry} and induction.
\end{proof}

\subsection{Estimating the cubic term}
Our goal here is to prove
\begin{equation}
\label{nonlin3}
\|h_2(r,u)\cdot v^3\|_{L^1_t\left(\dot{B}^1_{2,1}\cap\dot{B}^0_{2,1}(\R^{4})\right)}\,\lesssim \,\|v\|^3_{\tilde{X}},
\end{equation}
when
\begin{equation}
\label{tildex}
\|v\|_{\tilde{X}}\,:=\,\|\partial v\|_{L^{\infty}_t\left(\dot{B}^1_{2,1}\cap \dot{B}^{0}_{2,1}(\R^4)\right)}+\|\partial v\|_{L^2_t\left(\dot{B}^{1/6}_{6,1}\cap\dot{B}^{-5/6}_{6,1}(\R^4)\right)}
\end{equation}
is sufficiently small. 

Using  the Sobolev embedding \eqref{sob}, we deduce
\begin{equation}
\left\|\frac{\sin u}{r}\right\|_{L_{t,x}^{\infty}}\,\leq\, \|v\|_{L_{t,x}^{\infty}}\,\lesssim\, \|v\|_{L^\infty \dot{B}^2_{2,1}}\,\lesssim\, \|v\|_{\tilde{X}}\,\ll\, 1,
\label{phi}
\end{equation}
which allows us to infer
\begin{equation}
\label{cubser}
h_2(r,u)\cdot v^3\,=\,\sum\limits_{k=0}^{\infty}\,\tilde{h}_2(u)\,v^3\, (-1)^k\,\left( \frac{\sin u}{r}\right)^{2k}.
\end{equation}

First, we focus on the $L^1\,\dot{B}^1_{2,1}$ norm.
\begin{prop}
The following estimate is true for $\|v\|_{\tilde{X}}\,\ll\,1$:
\begin{equation}
\label{3b121}
\|h_2(r,u)\cdot v^3\|_{L^1\dot{B}^1_{2,1}}\,\lesssim\, \|v\|_{\ti X}^3.
\end{equation}
\end{prop}
\begin{proof}
If $k\geq 1$, the product bound \eqref{Bes-prod} implies
\begin{equation}
\label{3k}
\begin{aligned}
\Bigg\|\tilde{h}_2(u)\,v^3\, \left(\frac{\sin u}{r}\right)^{2k}\Bigg\|_{L^1\dot{B}^1_{2,1}}
\lesssim\, &\left\|\tilde{h}_2(u)\,v^3\right\|_{L^1\dot{B}^1_{2,1}}\,\left\|\left(\frac{\sin u}{r}\right)^{2k}\right\|_{L_{t,x}^{\infty}}\\&+\left\|\tilde{h}_2(u)\,v^3\right\|_{L^1L^{\infty}}\,\left\|\left(\frac{\sin u}{r}\right)^{2k}\right\|_{L^{\infty}\dot{B}^1_{2,1}}.
\end{aligned}
\end{equation}
Both norms involving $\frac{\sin u}{r}$ are estimated, based on \eqref{sob} and \eqref{ysin}, as
\begin{equation}
\label{linfsin}
\left\|\left(\frac{\sin u}{r}\right)^{2k}\right\|_{L_{t,x}^{\infty}}\,\lesssim\, \|v\|^{2k}_{L_{t,x}^{\infty}}\,\lesssim\, C^k\|v\|^{2k}_{\ti X}
\end{equation}
and
\begin{equation}
\label{b121sin}
\left\|\left(\frac{\sin u}{r}\right)^{2k}\right\|_{L^{\infty}\dot{B}^1_{2,1}}\,\lesssim\,\left\|\left(\frac{\sin u}{r}\right)^{2k}\right\|_{L^{\infty}Y} \,\lesssim\,C^k\|v\|^{2k}_{L^{\infty}Y}\,\lesssim\,C^k\|v\|^{2k}_{\ti X}.
\end{equation}

Next, using the decay estimates \eqref{analytic}, the embedding $\dot{B}^{2/3}_{6,1}(\R^4)\subset L^{\infty}(\R^4)$, and the interpolation relation
\begin{equation}
\left(\dot{B}^{7/6}_{6,1}, \dot{B}^{1/6}_{6,1}\right)_{\frac 12, 1}=\dot{B}^{2/3}_{6,1},
\end{equation}
we derive
\begin{equation}
\label{l1linfcub}
\begin{aligned}
\left\|\tilde{h}_2(u)v^3\right\|_{L^1L^{\infty}}\,&\lesssim\, \left\|v^3\right\|_{L^1L^{\infty}}\,\lesssim\, \|v\|_{L_{t,x}^{\infty}}\,\|v\|^2_{L^{2}L^{\infty}}\,\lesssim\, \|v\|_{\ti X}\,\|v\|^2_{L^{2}\dot{B}^{2/3}_{6,1}}\\&\lesssim\, \|v\|_{\ti X}\,\|v\|_{L^{2}\dot{B}^{7/6}_{6,1}}\,\|v\|_{L^{2}\dot{B}^{1/6}_{6,1}}\,\lesssim\, \|v\|^3_{\ti X}.
\end{aligned}
\end{equation}

Thus, we are left to investigate $\|\tilde{h}_2(u)\,v^3\|_{L^1\dot{B}^1_{2,1}}$, for which we rely on 
\begin{equation}
\left(\dot{H}^2, L^2\right)_{\frac{1}{2},1}\,=\,\dot{B}^1_{2,1}
\end{equation}
to infer
\begin{equation}
\label{v3l1b21}
\begin{aligned}
\|\tilde{h}_2(u)\, v^3\|_{L^1\dot{B}^1_{2,1}}\,\lesssim\, \|\tilde{h}_2(u)\, v^3\|_{L^1\dot{H}^2}^{1/2}\ \|\tilde{h}_2(u)\, v^3\|_{L^1L^2}^{1/2}.
\end{aligned}
\end{equation}

The second norm can be bounded, based on \eqref{analytic} and the Besov embeddings
\begin{equation}
\dot{B}^{1/6}_{6,1}(\R^4)\subset L^{8}(\R^4)\quad\text{and}\quad \dot{B}^{1}_{2,1}(\R^4)\subset L^4(\R^4),
\label{l8}
\end{equation}
as 
\begin{equation}
\label{v3l1l2}
\begin{aligned}
\|\tilde{h}_2(u)\, v^3\|_{L^1L^2}\,\lesssim\, \|v^3\|_{L^1L^2}\,&\lesssim\, \|v\|^2_{L^2L^8}\,\|v\|_{L^{\infty}L^4}\\&\lesssim\, \|v\|^2_{L^2\dot{B}^{1/6}_{6,1}}\,\|v\|_{L^{\infty}\dot{B}^{1}_{2,1}}\,\lesssim\,\|v\|^3_{\ti X}.
\end{aligned}
\end{equation}

Finally, taking advantage of \eqref{analytic}, we obtain
\EQ{
\left|\p_{rr}\left(\tilde{h}_2(u)\, v^3\right)\right|\,\lesssim\,\left|v^2\,v_{rr}\right|\,+\,\left|v\,v^2_{r}\right|\,+\,\left|v^3\,v_{r}\right|\,+\,\left|v^5\right|.}
Each of the terms on the right-hand side can be estimated in $L^1L^2$ through straightforward Besov embeddings as follows:
\begin{equation}
\label{v2vrr}
 \|v^2v_{rr}\|_{L^1L^2}\,\lesssim\,\|v\|^2_{L^2L^{\infty}}\,\|v_{rr}\|_{L^{\infty}L^2}\,\lesssim\,\|v\|^2_{\ti X}\,\|v\|_{L^{\infty}\dot{B}^2_{2,1}}\,\lesssim\,\|v\|^3_{\ti X},
\end{equation}
\begin{equation}
\label{v3vr}
\|v^3v_r\|_{L^1L^2}\,\lesssim\,\|v\|_{L^{\infty}L^{\infty}}\,\|v\|^2_{L^2L^{\infty}}\,\|v_r\|_{L^{\infty}L^2}\,\lesssim\,\|v\|^3_{\ti X}\,\|\p v\|_{L^{\infty}\dot{B}^0_{2,1}}\,\lesssim\,\|v\|^4_{\ti X},
\end{equation}
\begin{equation}
\label{vvr2}
\|vv_r^2\|_{L^1L^2}\,\lesssim\,\|v\|_{L^{\infty}L^4}\,\|v_r\|^2_{L^2L^8}\,\lesssim\,\|v\|_{\ti X}\,\|\p v\|^2_{L^2\dot{B}^{1/6}_{6,1}}\,\lesssim\,\|v\|^3_{\ti X},
\end{equation}
\begin{equation}
\label{v5}
\|v^5\|_{L^1L^2}\,\lesssim\,\|v\|^2_{L^\infty_{t,x}}\, \|v^3\|_{L^1L^2}\,\lesssim\,\|v\|^5_{\ti X}.
\end{equation}

As
\begin{equation}
\|\tilde{h}_2(u)\, v^3\|_{L^1\dot{H}^2}\,\simeq\, \|\p_{rr}\left(\tilde{h}_2(u)\, v^3\right)\|_{L^1L^2},
\end{equation}
the previous estimates imply
\begin{equation}
\|\tilde{h}_2(u)\, v^3\|_{L^1\dot{H}^2}\,\lesssim\,\|v\|^3_{\ti X},
\label{v3h2}
\end{equation} 
which, together with \eqref{v3l1b21} and \eqref{v3l1l2}, yields
\begin{equation}
\|\tilde{h}_2(u)\, v^3\|_{L^1\dot{B}^1_{2,1}}\,\lesssim\,\|v\|^3_{\ti X}.
\label{v3l1b121v2}
\end{equation} 
The argument is then concluded by combining the last estimate with \eqref{cubser}, \eqref{linfsin}, \eqref{b121sin}, and \eqref{l1linfcub}. 
\end{proof}

In order to finish the analysis of the cubic term, we need to estimate its $L^1\dot{B}^0_{2,1}$ norm.
\begin{prop}
The following estimate is true for $\|v\|_{\tilde{X}}\,\ll\,1$:
\begin{equation}
\label{3b021}
\|h_2(r,u)\cdot v^3\|_{L^1\dot{B}^0_{2,1}}\,\lesssim\, \|v\|_{\ti X}^3.
\end{equation}
\end{prop}
\begin{proof}
We rely on \eqref{b021} to infer, for $k\geq 1$,
\begin{equation}
\label{3k2}
\begin{aligned}
\Bigg\|\tilde{h}_2(u)\,v^3\, \left(\frac{\sin u}{r}\right)^{2k}\Bigg\|_{L^1\dot{B}^0_{2,1}}
\,\lesssim\, &\left\|\tilde{h}_2(u)\,v^3\right\|_{L^1\dot{B}^1_{2,1}}\,\left\|\left(\frac{\sin u}{r}\right)^{2k}\right\|_{L^{\infty}L^4}\\&+\,\left\|\tilde{h}_2(u)\,v^3\right\|_{L^1L^4}\,\left\|\left(\frac{\sin u}{r}\right)^{2k}\right\|_{L^{\infty}\dot{B}^1_{2,1}}.
\end{aligned}
\end{equation}

The first and the fourth norms on the right-hand side have already been estimated in the previous proposition.  Using again the Besov embeddings \eqref{l8}, we deduce
\begin{equation}
\aligned
\left\|\left(\frac{\sin u}{r}\right)^{2k}\right\|_{L^{\infty}L^4}\,&\lesssim\,\left\| \frac{\sin u}{r}\right \|^{2k-1}_{L_{t,x}^{\infty}}\,\left\|\frac{\sin u}{r}\right\|_{L^{\infty}L^4}\\
&\lesssim\,\left\|v\right\|^{2k-1}_{L_{t,x}^{\infty}}\,\left\|v\right\|_{L^{\infty}\dot{B}^1_{2,1}}\,\lesssim\,\|v\|_{\ti X}^{2k}
\endaligned
\label{sinl4}
\end{equation}
and 
\begin{equation}
\label{v3l4}
\left\|\tilde{h}_2(u)\,v^3\right\|_{L^1L^4}\,\lesssim\,\|v\|_{L_{t,x}^{\infty}}\,\|v\|^2_{L^2L^8}\,\lesssim\,\|v\|^3_{\ti X},
\end{equation}
which allow us to conclude \eqref{3b021}.
\end{proof}

\subsection{Estimating the quintic term}
For the quintic nonlinearity we are able to prove
\begin{equation}
\label{quintic}
\|h_3(r,u)\cdot v^5\|_{L^1\left(\dot{B}^1_{2,1}\cap\dot{B}^0_{2,1}(\R^{4})\right)}\,\lesssim\, \|v\|^5_{\ti X}.
\end{equation}
The argument has many similarities with the one for the cubic nonlinearity, as the extra $v^2$ present here in most expressions can be bounded directly in  $L^{\infty}_{t,x}$, which is controlled by $\|v\|_{\ti X}$. Moreover, the coefficients of the cubic and quintic terms, i.e., $\tilde{h}_2(u)$ and $\tilde{h}_3(u)$, verify almost identical decay estimates, according to \eqref{analytic}.

\subsection{Estimating the quartic term}
For this nonlinear term we intend to address first its $L^1\dot{B}^1_{2,1}$ norm. 
\begin{prop}
The following estimate is true for $\|v\|_{\tilde{X}}\,\ll\,1$:
\begin{equation}
\label{quartic}
\|h_1(r,u)\cdot v^3v_r\|_{L^1\dot{B}^1_{2,1}}\,\lesssim\, \|v\|^4_{\ti X}.
\end{equation}
\end{prop}
\begin{proof}
We proceed as before, based on \eqref{phi} and \eqref{Bes-prod}, to derive  
\begin{equation}
\label{qser}
h_1(r,u)\cdot v^3\,v_r\,=\,\sum\limits_{k=0}^{\infty}\,\tilde{h}_1(u)\,v^3\,v_r \,(-1)^k\,\left( \frac{\sin u}{r}\right)^{2k}
\end{equation}
and
\begin{equation}
\label{4k}
\begin{aligned}
\Bigg\|\tilde{h}_1(u)\,v^3v_r \left(\frac{\sin u}{r}\right)^{2k}\Bigg\|_{L^1\dot{B}^1_{2,1}}\,\lesssim\, &\left\|\tilde{h}_1(u)\,v^3v_r\right\|_{L^1\dot{B}^1_{2,1}}\,\left\|\left(\frac{\sin u}{r}\right)^{2k}\right\|_{L_{t,x}^{\infty}}\\
&+\,\left\|\tilde{h}_1(u)\,v^3v_r\right\|_{L^1L^4}\,\left\|\left(\frac{\sin u}{r}\right)^{2k}\right\|_{L^{\infty}\dot{B}^1_{4,1}}.
\end{aligned}
\end{equation}

In what concerns the norms for $\frac{\sin u}{r}$, the first one was investigated in \eqref{linfsin}, while the second one can be estimated using \eqref{Bes2} and \eqref{ysin} as
\begin{equation}
\aligned
\label{b141sin}
\left\|\left(\frac{\sin u}{r}\right)^{2k}\right\|_{L^{\infty}\dot{B}^1_{4,1}}\,&\lesssim\,\left\|\left(\frac{\sin u}{r}\right)^{2k}\right\|_{L^{\infty}\dot{B}^2_{2,1}}\,\lesssim\,\left\|\left(\frac{\sin u}{r}\right)^{2k}\right\|_{L^{\infty}Y}\\
&\lesssim\,C^k\|v\|^{2k}_{L^{\infty}Y}\,\lesssim\,C^k\|v\|^{2k}_{\ti X}.
\endaligned
\end{equation}
The Besov embedding \eqref{l8} allows us to control the $L^1L^4$ norm by 
\begin{equation}
\label{l1l4q}
\left\|\tilde{h}_1(u)\,v^3v_r\right\|_{L^1L^4}\,\lesssim\, \left\|v\right\|^2_{L_{t,x}^{\infty}}\,\|v\|_{L^{2}L^8}\,\|v_r\|_{L^{2}L^8}\,\lesssim\, \|v\|^4_{\ti X}.
\end{equation}

The analysis of the $L^1\dot{B}^1_{2,1}$ norm is the more intricate part of the proof for \eqref{quartic}, because we can not rely on the interpolation approach used for the cubic nonlinearity. The $\dot{H}^2_x$ norm appearing there would introduce a $v_{rrr}$ term for which we do not have good estimates. 

Instead, using \eqref{Bes-prod}, \eqref{l1linfcub}, and \eqref{l8}, we can infer
\begin{equation}
\aligned
\left\|\tilde{h}_1(u)\,v^3v_r \right\|_{L^1\dot{B}^1_{2,1}}\,&\lesssim\, \left\|\tilde{h}_1(u)\,v^3\right\|_{L^1L^\infty}\,\left\|v_r\right\|_{L^\infty\dot{B}^1_{2,1}}\,+\,\left\|\tilde{h}_1(u)\,v^3\right\|_{L^1\dot{B}^1_{4,1}}\,\left\|v_r\right\|_{L^{\infty}L^4}\\
&\lesssim\, \|v\|^4_{\ti X}\,+\,\left\|\tilde{h}_1(u)\,v^3\right\|_{L^1\dot{B}^1_{4,1}}\,\|v\|_{\ti X}.
\endaligned
\label{v3l1b21q}
\end{equation}

For the $L^1\dot{B}^1_{4,1}$ norm, we rely on the interpolation relations
\begin{equation}
\left(\dot{B}^{4/3}_{2,1}, \dot{B}^{8/9}_{6,1}\right)_{[\frac{3}{4}]}=\dot{B}^1_{4,1},\quad
\left(\dot{H}^2, L^2\right)_{\frac{1}{3},1}=\dot{B}^{4/3}_{2,1},\quad\ \left(\dot{W}^{1,6}, L^6\right)_{\frac{1}{9},1}=\dot{B}^{8/9}_{6,1},
\end{equation}
to deduce
\begin{equation}
\aligned
&\left\|\tilde{h}_1(u)\,v^3\right\|_{L^1\dot{B}^1_{4,1}}\lesssim\left\|\tilde{h}_1(u)\,v^3\right\|^{1/4}_{L^1\dot{B}^{4/3}_{2,1}}\,\left\|\tilde{h}_1(u)\,v^3\right\|^{3/4}_{L^1\dot{B}^{8/9}_{6,1}}\\
&\lesssim\left\|\tilde{h}_1(u)\,v^3\right\|^{1/6}_{L^1\dot{H}^2}\,\left\|\tilde{h}_1(u)\,v^3\right\|^{1/12}_{L^1L^2}\,\left\|\tilde{h}_1(u)\,v^3\right\|^{2/3}_{L^1\dot{W}^{1,6}}\,\left\|\tilde{h}_1(u)\,v^3\right\|^{1/12}_{L^1L^6}.
\label{qinterp}
\endaligned
\end{equation}

We claim that the $L^1\dot{H}^2$ and $L^1L^2$ norms have already been treated by \eqref{v3h2} and \eqref{v3l1l2}, respectively, due to similar decay properties for $\tilde{h}_1(u)$  and $\tilde{h}_2(u)$. Hence, we obtain 
\begin{equation}
\|\tilde{h}_1(u)\, v^3\|_{L^1\dot{H}^2}\,+\,\|\tilde{h}_1(u)\, v^3\|_{L^1L^2}\,\lesssim\,\|v\|^3_{\ti X}.
\label{qh2l2}
\end{equation} 
For the last norm, we use Besov embeddings and 
\begin{equation}
\left(\dot{B}^{7/6}_{6,1}, \dot{B}^{1/6}_{6,1}\right)_{\frac 23, 1}=\dot{B}^{1/2}_{6,1}
\end{equation}
to derive
\begin{equation}
\aligned
\|\tilde{h}_1(u)\, v^3\|_{L^1L^6}\,&\lesssim\,\|v\|_{L_{t,x}^{\infty}}\,\|v\|_{L^2L^8}\,\|v\|_{L^2L^{24}}\,\lesssim \,\|v\|^2_X\, \|v\|_{L^2\dot{B}^{1/2}_{6,1}}\\
&\lesssim\, \|v\|^2_{\ti X} \,\|v\|_{L^2\dot{B}^{7/6}_{6,1}}^{1/3}\, \|v\|_{L^2\dot{B}^{1/6}_{6,1}}^{2/3}\,\lesssim\, \|v\|_{\ti X}^3.
\endaligned
\label{ql1l6}
\end{equation}

Finally, in dealing with the $L^1\dot{W}^{1,6}$ norm, we rely on \eqref{analytic} to infer 
\begin{equation}
\left|\p_{r}\left(\tilde{h}_1(u)\, v^3\right)\right|\,\lesssim\,\left|v^2\,v_{r}\right| \,+\, \left|v^4\right|,
\end{equation}
while previous estimates imply
\begin{equation}
\|v^2\,v_r\|_{L^1L^6}\,\lesssim\,\|v\|_{L_{t,x}^{\infty}}\,\|v\|_{L^2L^{24}}\,\|v_r\|_{L^2L^8}\,\lesssim \,\|v\|^3_{\ti X}
\end{equation}
and
\begin{equation}
\|v^4\|_{L^1L^6}\,\lesssim\,\|v\|_{L_{t,x}^{\infty}}\,\|v^3\|_{L^1L^{6}}\,\lesssim \,\|v\|^4_{\ti X}.
\end{equation}
As
\begin{equation}
\|\tilde{h}_1(u)\, v^3\|_{L^1\dot{W}^{1,6}}\,\simeq\, \|\p_{r}\left(\tilde{h}_1(u)\, v^3\right)\|_{L^1L^6},
\end{equation}
we deduce
\begin{equation}
\|\tilde{h}_1(u)\, v^3\|_{L^1\dot{W}^{1,6}}\,\lesssim\, \|v\|^3_{\ti X},
\label{qw16}
\end{equation}
for $\|v\|_X$ sufficiently small. Therefore, based on the previous bound, \eqref{v3l1b21q}, \eqref{qinterp}, \eqref{qh2l2}, and \eqref{ql1l6}, we obtain
\EQ{\left\|\tilde{h}_1(u)\,v^3v_r \right\|_{L^1\dot{B}^1_{2,1}}\,\lesssim\,  \|v\|^4_{\ti X},\label{v3vrtih}} 
which concludes the proof of \eqref{quartic}.
\end{proof}

The last ingredient needed in finishing the proof of \eqref{nonlin*} is a favorable estimate for the $L^1\dot{B}^0_{2,1}$ norm of the quartic nonlinearity.
\begin{prop}
The following estimate is true for $\|v\|_{\tilde{X}}\,\ll\,1$:
\begin{equation}
\label{4b021}
\|h_1(r,u)\cdot v^3v_r\|_{L^1\dot{B}^0_{2,1}}\,\lesssim\, \|v\|_{\ti X}^4.
\end{equation}
\end{prop}
\begin{proof}
As in the analysis of the corresponding norm for the cubic term, we use \eqref{b021} to infer, for $k\geq 1$, \begin{equation}
\label{4k2}
\begin{aligned}
\Bigg\|\tilde{h}_1(u)\,v^3v_r\, \left(\frac{\sin u}{r}\right)^{2k}\Bigg\|_{L^1\dot{B}^0_{2,1}}
\,\lesssim\, &\left\|\tilde{h}_1(u)\,v^3v_r\right\|_{L^1\dot{B}^1_{2,1}}\,\left\|\left(\frac{\sin u}{r}\right)^{2k}\right\|_{L^{\infty}L^4}\\&+\,\left\|\tilde{h}_1(u)\,v^3v_r\right\|_{L^1L^4}\,\left\|\left(\frac{\sin u}{r}\right)^{2k}\right\|_{L^{\infty}\dot{B}^1_{2,1}}.
\end{aligned}
\end{equation}
There is nothing new left to argue in proving \eqref{4b021}, because all of the four norms on the right-hand side have been previously bounded: the first by \eqref{v3vrtih}, the second by \eqref{sinl4}, the third by \eqref{l1l4q}, and the final one by \eqref{b121sin}. 
\end{proof}

\section{Estimating the null-form nonlinearity}

In this section we complete the proof for our main result, Theorem \ref{mainth}, as follows. 

First, we define the solution space $X$ and, since we are at the critical regularity level for $4+1$-dimensional wave maps, it needs to factor in the null structure of the last nonlinearity. This is why we cannot use $\ti X$, which was defined in the previous section by \eqref{tildex}, as an iteration space. Using the properties of $X$, we prove the homogeneous estimates \eqref{hom}, 
\begin{equation}
\|\partial v\|_{L^{\infty}_t\left(\dot{B}^1_{2,1}\cap \dot{B}^{0}_{2,1}(\R^4)\right)}\,+\,\|\partial v\|_{L^2_t\left(\dot{B}^{1/6}_{6,1}\cap\dot{B}^{-5/6}_{6,1}(\R^4)\right)}\,\lesssim\,\|v\|_X, 
\label{xdom}
\end{equation}
and
\begin{equation}
\|\Box^{-1}H\|_{X}\,\lesssim\,\|H\|_{L^1_t\left(\dot{B}^1_{2,1}\cap\dot{B}^0_{2,1}(\R^{4})\right)}.
\label{BoxXmin}
\end{equation}
Together with \eqref{nonlin*}, the last two inequalities imply \eqref{nonlin3*}, which is the desired bound for ${\ti N}(v)$.

Secondly, we demonstrate the trilinear estimate 
\begin{equation}
\|r^2\,u\,v\,w\|_X\,\lesssim\, \|u\|_X\, \|v\|_X\, \|w\|_X, \label{tril}
\end{equation}
which is the main new contribution of this article. This allows us to bound the null-form nonlinearity \eqref{nf} by
\begin{equation}
\label{nonlinnf}
\|\Box^{-1}\left(h_4(r,u)\cdot v(v_r^2-v_t^2)\right)\|_{X}\,\lesssim\, \|v\|^3_X,
\end{equation}
when $\|v\|_X$ is sufficiently small, and finish the argument.

\subsection{The solution space $X$ and its properties} 
Using $\chi\in C_0^\I(\R)$, which is a smooth cutoff function satisfying
\EQ{
 \supp\chi\subset(1/2,2), \pq \pq \sum_{\la\in 2^\Z}\chi(\la^{-1}s)=1, \ (\forall) s\neq 0,}
we define the spacetime Fourier multipliers
\begin{equation}
\aligned
 &A_\la(D)=\F^{-1}\,\chi\left(\la^{-1}|(\tau,\xi)|\right)\,\F,\\ 
 B_\la(D)=\F^{-1}\,\chi&\left(\la^{-1}\frac{|\tau^2-|\xi|^2|}{|(\tau,\xi)|}\right)\,\F,\qquad
\tilde{B}_\mu(D)=\sum_{j\geq 0} B_{2^{-j}\mu}(D),
\endaligned
\end{equation}
where $\la$, $\mu \in 2^\mathbb{Z}$, $\F$ is the Fourier transform in $(t,x)\in \R^{n+1}$, and $|(\tau,\xi)|=\sqrt{\tau^2+|\xi|^2}$. Next, for functions $w=w(t,x)$ whose Fourier transform is supported at frequency $|(\tau,\xi)|\simeq \la$, we can associate the norms
\begin{equation}
\|w\|_{X_\la^s}=\sum_{\mu\in 2^\mathbb{Z}} \mu^s\,\|B_\mu(D)w\|_{L^2_{t,x}}, \qquad \|w\|_{Y_\la}=\|w\|_{L^\infty L^2} + \la^{-1}\|\Box w\|_{L^1L^2}.
\end{equation}
 
In \cite{T}, Tataru introduced the Besov-type hyperbolic spaces $F$ and $\Box F$, which are described by 
\begin{equation}
\label{F}
\|w\|_F\, =\,\sum\limits_{\la\in 2^\mathbb{Z}}\,\la^{n/2}\,\|A_{\la}(D)w\|_{F_{\la}}, \qquad F_{\la}\,=\,X^{1/2}_{\la}\,+\,Y_{\la},
\end{equation}
and 
\begin{equation}
\label{Box F}
\|w\|_{\Box F}\, =\,\sum\limits_{\la\in 2^\mathbb{Z}}\,\la^{n/2}\,\|A_{\la}(D)w\|_{\Box F_{\la}}, \qquad
\Box F_{\la}\,=\,\la \left(X^{-1/2}_{\la}\,+\,(L^1L^2)_\la\right),
\end{equation}
respectively. We collect in the next proposition the properties proved there for these function spaces.
\begin{prop}(\cite{T})\ Let $n\geq 4$. 

i) For the linear wave equation, the following estimates are true:
\EQ{
\|S(v_0,v_1)\|_F\,\lesssim\,\|(v_0,v_1)\|_{\dot{B}^{n/2}_{2,1}\times \dot{B}^{(n-2)/2}_{2,1}},\label{S}}
\EQ{
\|\Box^{-1}H\|_{F}\,\lesssim\,\|H\|_{\Box F}\,\lesssim\,\|H\|_{L^1\dot{B}^{(n-2)/2}_{2,1}}.
\label{Box1}}

ii) If $(q,r)$ is a Strichartz wave-admissible pair, i.e.,
\begin{equation}
2\,\leq \,q,\,r\,\leq \,\infty \qquad \text{and} \qquad \frac{2}{q}+\frac{n-1}{r}\leq \frac{n-1}{2},
\end{equation}
then 
\begin{equation}
\la^{n/r+ 1/q-n/2}\,\|v\|_{L^q L^r}\,\lesssim\,\|v\|_{F_\la}
\label{str}
\end{equation}
holds uniformly in $\la$.

iii) The following bilinear inequalities are true:
\begin{equation}
\|v\, w\|_{F}\,\lesssim\, \|v\|_{F}\, \|w\|_{F}, \qquad
\|v\, w\|_{\Box F}\,\lesssim\, \|v\|_{F}\, \|w\|_{\Box F}.  \label{algF}
\end{equation}
As a consequence,
\begin{equation}
\|v_t\,w_t\,-\,\nabla_xv\,\nabla_xw\|_{\Box F}\,\lesssim\, \|v\|_{F}\, \|w\|_{F}. 
\label{nullF1}
\end{equation}

iv) If $ w_\la=A_\la(D)w$ and
\EQ{w_\la=w_\la^{<\mu}+w_\la^{>\mu}, \qquad w_\la^{<\mu} = \tilde{B}_\mu(D)w_\la,\label{wlm}}
then
\begin{equation}
\|w_\la^{<\mu}\|_{F_\la} \lesssim \|w_\la\|_{F_\la}, \quad \|w_\la^{<\mu}\|_{\Box F_\la} \lesssim \|w_\la\|_{\square F_\la}, \quad \|w_\la^{>\mu}\|_{L^1L^2} \lesssim \mu^{-1} \|w_\la\|_{Y_\la}, \label{tiB}
\end{equation}
hold uniformly for $\mu\leq \la$.
\end{prop}

\begin{rem}
For the radial linear wave equation, Sterbenz \cite{MR2128434} was able to enlarge the set of Strichartz-admissible pairs to include the ones which satisfy
\begin{equation}
\frac{1}{q}+\frac{n-1}{r} < \frac{n-1}{2}.
\label{improvstr}
\end{equation}
As the argument for the dyadic bounds \eqref{str} is based on a straightforward application of the classical Strichartz inequalities, we could modify it to prove that \eqref{str} hold in the extended range \eqref{improvstr} for radial functions.
\end{rem}

This was followed by Geba-Nakanishi-Rajeev \cite{GNR} who introduced the function spaces $|\nabla| F$ and $|\nabla| \Box F$, which are given by 
\begin{equation}
\label{nabF}
\|w\|_{|\nabla| F}\, =\,\sum\limits_{\la\in 2^\mathbb{Z}}\,\la^{(n-2)/2}\,\|A_{\la}(D)w\|_{F_{\la}}, 
\end{equation}
and
\begin{equation}
\label{BoxgradF}
\|w\|_{|\nabla| \Box F}\,=\,\sum\limits_{\la\in 2^\mathbb{Z}}\,\la^{(n-2)/2}\,\|A_{\la}(D)w\|_{\Box F_{\la}},
\end{equation}
respectively. In that paper, they used also the radial space
\begin{equation}
\label{Zdef}
 Z\,=\,\left\{v \in \mathcal{S}'(\R^{n+1})\big|\, v=v(t,x)=v(t,r),\ \|v\|_{F\cap|\nabla|F}<\infty \right\}. 
\end{equation}
For these spaces, arguments similar to the ones in \cite{T} lead to:
\begin{prop}(\cite{GNR})\ For $n\ge 4$,
\begin{equation}
\|S(v_0,v_1)\|_{|\nabla| F}\,\lesssim\,\|(v_0,v_1)\|_{\dot{B}^{(n-2)/2}_{2,1}\times \dot{B}^{(n-4)/2}_{2,1}},
\label{S2}
\end{equation}
\begin{equation}
\|\Box^{-1}H\|_{|\nabla| F}\,\lesssim\,\|H\|_{|\nabla|\Box F}\,\lesssim\,\|H\|_{L^1\dot{B}^{(n-4)/2}_{2,1}},
\label{Box2}
\end{equation}
\begin{equation}
\|v\, w\|_{|\nabla| F}\,\lesssim\, \|v\|_{F}\, \|w\|_{|\nabla| F}, \qquad \|v\,w\|_{|\nabla|\Box F}\,\lesssim\, \|v\|_{|\nabla| F}\, \|w\|_{\Box F} \label{algnbF},
\end{equation}
\begin{equation}
\|v_t\,w_t\,-\,\nabla_xv\,\nabla_xw\|_{|\nabla|\Box F}\,\lesssim\, \|v\|_{F\cap |\nabla| F}\, \|w\|_{F\cap |\nabla| F}, 
\label{nullF}
\end{equation}
hold, while if $n\ge 5$ and $v=v(t,r)$ and $w=w(t,r)$ are two radial functions, then
\EQ{ \label{Z5}
 \|r \,v\, w\|_{Z}\, \lec\, \|v\|_{Z}\, \|w\|_{Z}.}
\end{prop}

We will use as a solution space the $n=4$ version of $Z$, 
\begin{equation}
\label{Xdef}
 X:\,=\,\left\{v \in \mathcal{S}'(\R^{4+1})\big|\, v=v(t,x)=v(t,r),\ \|v\|_{F\cap|\nabla|F}<\infty \right\}. 
\end{equation} 

First, due to \eqref{algF} and \eqref{algnbF}, we deduce that $X$ is an algebra, i.e.,
\begin{equation}
\|v\, w\|_{X}\,\lesssim\, \|v\|_{X}\, \|w\|_{X}.
\label{algX}
\end{equation} 

Next,  \eqref{S} and \eqref{S2} imply together
\begin{equation}
\|S(v_0,v_1)\|_X\,\lesssim\,\|v_0\|_{\dot{B}^{2}_{2,1}\cap\dot{B}^{1}_{2,1} (\R^4)}+\|v_1\|_{\dot{B}^{1}_{2,1}\cap\dot{B}^{0}_{2,1}(\R^4)}\,=\,\|(v_0,v_1)\|_D,
\label{hom2}
\end{equation}
which is one of the homogeneous estimates \eqref{hom}. Also, combining \eqref{Box1} and \eqref{Box2}, we obtain \eqref{BoxXmin}.

As $(\infty, 2)$ and $(2,6)$ are Strichartz-admissible pairs for $n=4$, we derive, based on the definitions of $F$ and $|\nabla| F$,
\begin{equation}
\|\partial v\|_{L^{\infty}_t\left(\dot{B}^1_{2,1}\cap \dot{B}^{0}_{2,1}(\R^4)\right)}\,+\,\|\partial v\|_{L^2_t\left(\dot{B}^{1/6}_{6,1}\cap\dot{B}^{-5/6}_{6,1}(\R^4)\right)}\,\lesssim\,\|v\|_{F\cap |\nabla| F} 
\label{ll}.
\end{equation} 
In the radial case, this is precisely \eqref{xdom}. Moreover, the other half of  \eqref{hom},
\begin{equation}
\|w\|_{L_t^\infty\,D}\,=\,\|w\|_{L^\infty(\dot{B}^{2}_{2,1}\cap\dot{B}^{1}_{2,1} (\R^4))}+\|\p_t w\|_{L^\infty(\dot{B}^{1}_{2,1}\cap\dot{B}^{0}_{2,1}(\R^4))}\,\lesssim\,\|w\|_X,
\label{embed2}
\end{equation}
is an immediate consequence.

In what concerns the null-form nonlinearity, coupling \eqref{nullF1} and \eqref{nullF} leads to 
\begin{equation}
\|v_t^2-v_r^2\|_{\Box F \cap |\nabla| \Box F}\,\lesssim\, \|v\|^2_{Z}. 
\label{nullZ}
\end{equation}

\subsection{The trilinear estimate and the conclusion of the argument}
\begin{thm}
The following trilinear inequality holds:
\begin{equation}
\|r^2\,u\,v\,w\|_X\,\lesssim\, \|u\|_X\, \|v\|_X\, \|w\|_X.
\label{tril2}
\end{equation}
\end{thm}
\begin{proof}
For $n=4$, $q=\I$, and $\al=p$, the radial Sobolev inequality \eqref{pq} becomes
\begin{equation}
 \|r\phi_\mu\|_{L^\I(\R^4)} \lesssim \mu^{(4-p)/p} \|\phi\|_{L^p(\R^4)}.
\end{equation}
Together with \eqref{str}  in the improved radial range \eqref{improvstr}, it leads to
\begin{equation}
\label{rv}
\|r\,v\|_{L^q L^\infty}\, \lesssim\, \la^{1-1/q}\,\|v\|_{F_\la}, \qquad (\forall)\,2 < q \leq \infty.
\end{equation}

Given that $X$ is a Besov-type space with $l^1$ summability over dyadic decompositions, it is sufficient to prove \eqref{tril2} for single dyadic pieces, i.e.,
\EQ{u_\nu=A_\nu(D)u, \quad v_\mu=A_\mu(D)v, \quad w_\la=A_\la(D)w,  \qquad \nu, \,\mu, \,\la \in 2^\mathbb{Z},}
where we can assume, by symmetry, that $ \nu\leq \mu\leq \la$. We further decompose $w_\la$ as  
\begin{equation}
w_\la=w_\la^{<\mu}+w_\la^{>\mu}, \qquad w_\la^{<\mu} = \tilde{B}_\mu(D)w_\la,
\end{equation} 
and rely on \eqref{rv} and \eqref{tiB} to estimate $r^2 u_\nu v_\mu w_\la^{<\mu}$ as follows:
\begin{equation}
 \aligned
 \|r^2 \,u_\nu \,v_\mu \,w_\la^{<\mu}\|_{L^2_{t,x}}
 &\lesssim\, \|r \,u_\nu\|_{L^q L^\infty}\, \|r\, v_\mu\|_{L^{2q/(q-2)} L^\infty}\, \|w_\la^{<\mu}\|_{L^\infty L^2}
 \\
 &\lesssim\, \nu^{1-1/q} \,\mu^{1/2+1/q}\,\|u_\nu\|_{F_\nu}\,\|v_\mu\|_{F_\mu}\, \|w_\la\|_{F_\la},
\endaligned
\end{equation}
for $2<q<\infty$. We notice here the flexibility we have in choosing $q$ to satisfy
\[
\nu^{1-1/q}\, \leq\, 2\nu,\]
e.g.,
\[
q\,=\,\max\left\{2,\log_2 \left(\nu^{-1}\right)\right\}\,+\,1.
\]
As $\mathcal{F}\left(r^2 u_\nu v_\mu w_\la^{<\mu}\right)$ is supported in the set
\[
 |\tau|+|\xi|\, \lesssim\, \la, \qquad  ||\tau|-|\xi||\, \lesssim\, \mu,\]
we deduce
\begin{equation}
\label{x1}
 \|r^2\, u_\nu \,v_\mu \,w_\la^{<\mu}\|_{X^{1/2}_\la}\, \lesssim\, \nu \,(\mu^{3/2}+\mu)\, \|u_\nu\|_{F_\nu}\, \|v_\mu\|_{F_\mu}\, \|w_\la\|_{F_\la}.
\end{equation}

The analysis for  $r^2 \,u_\nu\, v_\mu \,w_\la^{>\mu}$ is clearly relevant only if $\mu\ll\la$. We need the decomposition 
\begin{equation}
w_\la^{>\mu} = \sum_{2^{4}\mu\,<\,d\,\leq\, \la} w_\la^d \,+ \,w_\la^0,
\end{equation}
where the spacetime Fourier support of $w_\la^d$ lies in $|\tau^2-|\xi|^2|\sim d\la$ and 
\begin{equation}
 \|w_\la^{>\mu}\|_{F_\la} \,\sim \,\sum_d\, \|w_\la^d\|_{X^{1/2}_\la} \,+ \,\|w_\la^0\|_{Y_\la}.
\end{equation}
The $X^{1/2}_\la$ component is estimated directly, using again \eqref{rv} and \eqref{tiB}:
 \begin{equation}
 \label{x2}
 \aligned
 \|r^2\, u_\nu \,v_\mu \,\sum_d\, w_\la^d\|_{X^{1/2}_\la} 
 &\lesssim \,\sum_d\, d^{1/2}\,\|r^2\, u_\nu\, v_\mu\, w_\la^d\|_{L^2_{t,x}}\\ &\lesssim \,\|r \,u_\nu\|_{L^\infty_{t,x}} \,\|r\, v_\mu\|_{L^\infty_{t,x}}\, \sum_d\, d^{1/2}\,\|w_\la^d\|_{L^2_{t,x}}\\ 
 &\lesssim\, \nu\, \mu \,\|u_\nu\|_{F_\nu}\, \|v_\mu\|_{F_\mu} \,\|w_\la\|_{F_\la}.
 \endaligned
 \end{equation}
For the $Y_\la$ component, a similar argument allows us to infer first that 
\begin{equation}
\label{x3}
\aligned
\|r^2 \,u_\nu \,v_\mu \,w_\la^0\|_{L^\infty L^2}\,
 &\lesssim\, \|r \,u_\nu\|_{L^\infty_{t,x}}\, \|r \,v_\mu\|_{L^\infty_{t,x}} \,\|w_\la^0\|_{L^\infty L^2}\\
 &\lesssim\, \nu\, \mu \,\|u_\nu\|_{F_\nu}\, \|v_\mu\|_{F_\mu}\,\| w_\la\|_{F_\la}.
\endaligned
\end{equation}
Secondly, as \eqref{tiB} implies
\begin{equation}
\label{w0}
 \|w_\la^0\|_{L^1L^2} \,\lesssim\, \mu^{-1}\|w_\la^0\|_{Y_\la},
\end{equation} 
it follows that
\begin{equation}
\label{x4}
\aligned
\la^{-1}\|\Box(r^2 \,u_\nu \,&v_\mu \,w_\la^0)\|_{L^1L^2}\\
&\lesssim  \nu\, \mu \,\|u_\nu\|_{F_\nu} \,\|v_\mu\|_{F_\mu}\,\| w_\la\|_{F_\la} 
\,+\, \la^{-1}\,\|\left[\Box, r^2\, u_\nu\, v_\mu\right] w_\la^0\|_{L^1L^2}\\
&\lesssim\,  \nu\, \mu \,\|u_\nu\|_{F_\nu} \,\|v_\mu\|_{F_\mu}\,\| w_\la\|_{F_\la} 
\,+\, \la^{-1}\, \mu\,\la \,\|r^2\,u_\nu\, v_\mu\|_{L^\infty_{t,x}}\, \|w_\la^0\|_{L^1L^2}\\
&\lesssim\, \nu\, \mu \,\|u_\nu\|_{F_\nu} \,\|v_\mu\|_{F_\mu}\,\| w_\la\|_{F_\la}.
 \endaligned
 \end{equation}

Connecting now \eqref{x1}, \eqref{x2}, \eqref{x3}, and \eqref{x4}, we obtain
\begin{equation}
\|r^2 u_\nu v_\mu w_\la\|_X\, \lesssim \,\nu  \,(\mu^2+\mu)\,(\la^2+\la) \|u_\nu\|_{F_\nu}\, \|v_\mu\|_{F_\mu}\, \|w_\la\|_{F_\la},
\end{equation}
which yields \eqref{tril2} by summing up over the three indices. 
\end{proof}

\begin{rem}
An interesting open problem is whether the bilinear estimate \eqref{Z5} is true when $n=4$. If that is the case, then it would imply \eqref{tril2}. The argument for proving \eqref{Z5} in \cite{GNR} can not be applied when $n=4$, because the $q=2$ endpoint for \eqref{rv} is unavailable. 
\end{rem}

By induction, the trilinear estimate \eqref{tril2} yields, for all $j\geq 0$,
\begin{equation}
\label{rX}
\|u_1\, u_2\, \ldots\, u_{2j+1}\|_{rX}\, \lesssim\, C^j\, \|u_1\|_{rX} \, \|u_2\|_{rX} \ldots \|u_{2j+1}\|_{rX},
\end{equation}
where $\|w\|_{rX}=\|w/r\|_X$. Hence, for any fixed $\alpha\in\mathbb{R}$,
\begin{equation}
\aligned
\left\|\frac{\sin \alpha u}{r}\right\|_X\,&\lesssim \,\sum_{j=0}^\infty \frac{|\alpha|^{2j+1}}{{(2j+1)!}}\,\left\|\frac{u^{2j+1}}{r}\right\|_X\,= \,\sum_{j=0}^\infty \frac{|\alpha|^{2j+1}}{{(2j+1)!}}\,\left\|u^{2j+1}\right\|_{rX}\\&\lesssim \,\sum_{j=0}^\infty \frac{C^{j}\,|\alpha|^{2j+1}}{{(2j+1)!}}\,\left\|u\right\|^{2j+1}_{rX}\,=\,\sum_{j=0}^\infty \frac{C^{j}\,|\alpha|^{2j+1}}{{(2j+1)!}}\,\left\|v\right\|^{2j+1}_{X}\,\lesssim\,\|v\|_X,
\endaligned
\label{sinu}
\end{equation}
for $\|v\|_X$ sufficiently small. As $X$ is an algebra and
\begin{equation}
\label{exph4}
 h_4(r,u)\cdot v\,=\,\frac{\sin 2u}{2r}\,\cdot\,\sum_{k=0}^\infty (-1)^k\,\left(\frac{\sin u}{r}\right)^{2k},
 \end{equation}
the previous inequality implies
\begin{equation}
\label{h4X}
\|h_4(r,u)\cdot v\|_X\,\lesssim\,\left\|\frac{\sin 2u}{r}\right\|_X\,\sum_{k=0}^\infty \left(C\left\|\frac{\sin u}{r}\right\|_X\right)^{2k}\,\lesssim\,\|v\|_X.
 \end{equation}
Together with \eqref{algF}, \eqref{algnbF}, and \eqref{nullZ}, it proves
\begin{equation}
\label{nullv2}
\|h_4(r,u)\cdot v(v_r^2-v_t^2)\|_{\Box F \cap |\nabla| \Box F}\,\lesssim\,\|v\|^3_X,
\end{equation}
which, based on \eqref{Box1} and \eqref{Box2}, provides us with the null-form estimate \eqref{nonlinnf}, thus finishing the argument.

\section*{Acknowledgements}
The first and the third author were supported in part by the National Science Foundation Career grant DMS-0747656. 

\bibliographystyle{amsplain}
\bibliography{anwb}

\providecommand{\bysame}{\leavevmode\hbox to3em{\hrulefill}\thinspace}
\providecommand{\MR}{\relax\ifhmode\unskip\space\fi MR }
\providecommand{\MRhref}[2]{%
  \href{http://www.ams.org/mathscinet-getitem?mr=#1}{#2}
}
\providecommand{\href}[2]{#2}
\begin{thebibliography}{10}

\bibitem{BS98}
R.~A. Battye and P.~M. Sutcliffe, \emph{Knots as stable soliton solutions in a
  three-dimensional classical field theory}, Phys. Rev. Lett. \textbf{81}
  (1998), no.~22, 4798--4801.

\bibitem{BS99}
\bysame, \emph{Solitons, links and knots}, R. Soc. Lond. Proc. Ser. A Math.
  Phys. Eng. Sci. \textbf{455} (1999), no.~1992, 4305--4331.

\bibitem{F75}
L.~D. Faddeev, \emph{Quantization of solitons}, Preprint IAS print-75-QS70
  (Inst. Advanced Study, Princeton, NJ, 1975), 32 pp.

\bibitem{F76}
\bysame, \emph{Some comments on the many-dimensional solitons}, Lett. Math.
  Phys. \textbf{1} (1976), no.~4, 289--293.

\bibitem{FN97}
L.~D. Faddeev and A.~J. Niemi, \emph{Stable knot-like structures in classical
  field theory}, Nature \textbf{387} (1997), 58--61.

\bibitem{GNR}
D.-A. Geba, K.~Nakanishi, and S.~G. Rajeev, \emph{Global well-posedness and
  scattering for {S}kyrme wave maps}, Commun. Pure Appl. Anal. \textbf{11}
  (2012), no.~5, 1923--1933.

\bibitem{LLZ}
Z.~Lei, F.~Lin, and Y.~Zhou, \emph{Global solutions of the evolutionary
  {F}addeev model with small initial data}, Acta Math. Sin. (Engl. Ser.)
  \textbf{27} (2011), no.~2, 309--328.

\bibitem{LY04}
F.~Lin and Y.~Yang, \emph{Existence of energy minimizers as stable knotted
  solitons in the {F}addeev model}, Comm. Math. Phys. \textbf{249} (2004),
  no.~2, 273--303.

\bibitem{LY104}
\bysame, \emph{The {F}addeev knots as stable solitons: existence theorems},
  Sci. China Ser. A \textbf{47} (2004), no.~2, 187--197.

\bibitem{S1}
T.~H.~R. Skyrme, \emph{A non-linear field theory}, Proc. Roy. Soc. London Ser.
  A \textbf{260} (1961), 127--138.

\bibitem{S2}
\bysame, \emph{Particle states of a quantized meson field}, Proc. Roy. Soc.
  Ser. A \textbf{262} (1961), 237--245.

\bibitem{S3}
\bysame, \emph{A unified field theory of mesons and baryons}, Nuclear Phys.
  \textbf{31} (1962), 556--569.

\bibitem{MR2128434}
J.~Sterbenz, \emph{Angular regularity and {S}trichartz estimates for the wave
  equation}, Int. Math. Res. Not. (2005), no.~4, 187--231, With an appendix by
  I. Rodnianski.

\bibitem{T}
D.~Tataru, \emph{Local and global results for wave maps. {I}}, Comm. Partial
  Differential Equations \textbf{23} (1998), no.~9-10, 1781--1793.

\bibitem{VK79}
A.~F. Vakulenko and L.~V. Kapitanski, \emph{Stability of solitons in {$S^{2}$}
  of a nonlinear {$\sigma $}-model}, Dokl. Akad. Nauk SSSR \textbf{246} (1979),
  no.~4, 840--842.

\end{thebibliography}

\end{document}